\documentclass{article}
\usepackage{amssymb}
\usepackage{verbatim}

\title{Models in which every nonmeager set is nonmeager in a nowhere
dense Cantor set}

\author{
{\sc Maxim R.~Burke} \thanks{Research supported by NSERC. The
author thanks the Department of Mathematics at the University of
Wisconsin for its hospitality during the academic year 1996/1997
when the earlier result mentioned in the introduction was produced
and the Department of Mathematics at the University of Toronto for
its hospitality during the academic year 2003/2004 when the
present paper was completed.}
\\
{\footnotesize Department of Mathematics and Statistics} \\
{\footnotesize University of Prince Edward Island} \\
{\footnotesize Charlottetown PE, Canada C1A 4P3} \\
{\footnotesize burke@upei.ca}
\and {\sc Arnold W. Miller} \thanks{Thanks to the Fields Institute
for Research in Mathematical Sciences for their support during
part of the time this paper was written and to Juris Steprans who
directed the special program in set theory and analysis.
\endgraf Key words and phrases: Property of Baire, Lebesgue measure,
Cantor set, oracle forcing
}
\\
{\footnotesize University of Wisconsin-Madison} \\
{\footnotesize Department of Mathematics, Van Vleck Hall} \\
{\footnotesize 480 Lincoln Drive} \\
{\footnotesize Madison, Wisconsin 53706-1388} \\
{\footnotesize miller@math.wisc.edu}\\
{\footnotesize http://www.math.wisc.edu/$\sim$miller}
}

\date{\today}

\def\a{\alpha}
\def\barM{\overline{M}}

\def\cof{\mathop{\rm cof}}
\def\d{\delta}
\def\dom{\mathop{\rm dom}}

\def\g{\gamma}

\def\la{\langle}
\def\m{\medskip}
\def\n{\noindent}
\def\o{\omega}

\def\Q{{\mathbb Q}}
\def\ra{\rangle}
\def\range{\mathop{\rm ran}}
\def\restr{\upharpoonright}
\def\R{{\mathbb R}}
\def\s{\sigma}

\def\sm{\setminus}
\def\sq{\subseteq}
\def\t{\theta}
\def\trap{\mathop{\rm Trap}}
\def\Trap{\trap(\barM)}

\def\proof{\noindent {\sc Proof. }}
\def\poset{{\mathbb P}}
\def\res{\upharpoonright}
\def\qed{\hspace{\stretch{1}}$\square$}

\def\perf{\mathop{\rm perfect}}
\def\irr{\mathop{\rm baire}}
\def\cantor{\mathop{\rm cantor}}

\newtheorem{theorem}{Theorem}[section]
\newtheorem{lemma}[theorem]{Lemma}
\newtheorem{proposition}[theorem]{Proposition}
\newtheorem{fact}[theorem]{\sc Fact}
\newtheorem{corollary}[theorem]{Corollary}
\newtheorem{claim}[theorem]{Claim}

\newtheorem{problem}[theorem]{Problem}
\newtheorem{definition}[theorem]{Definition}
\newtheorem{remark}[theorem]{Remark}

\newcommand{\thm}[2]{\begin{theorem}\label{#1}{\sl #2}\end{theorem}}

\newcommand{\prop}[2]{\begin{proposition}\label{#1}{\sl
#2}\end{proposition}}
\newcommand{\lem}[2]{\begin{lemma}\label{#1}{\sl #2}\end{lemma}}
\newcommand{\clm}[2]{\begin{claim}\label{#1}{\rm #2}\end{claim}}
\newcommand{\pr}[2]{\begin{problem}\label{#1}{\rm #2}\end{problem}}
\newcommand{\faact}[2]{\begin{fact}\label{#1}{\rm #2}\end{fact}}

\newcommand{\defi}[2]{\begin{definition}\label{#1}{\rm
#2}\end{definition}}
\newcommand{\rem}[2]{\begin{remark}\label{#1}{\rm #2}\end{remark}}

\begin{document}

\maketitle


\begin{center}
Abstract
\end{center}
\begin{quote}
We prove that it is relatively consistent with ZFC that in any
perfect Polish space, for every nonmeager set $A$ there exists a
nowhere dense Cantor set $C$ such that $A\cap C$ is nonmeager in
$C$. We also examine variants of this result and establish a
measure theoretic analog.
\end{quote}

\section{Introduction}

Our starting point is the following question of Laczkovich:
\begin{quote}
Does there exist (in ZFC) a nonmeager set that is relatively
meager in every nowhere dense perfect set?
\end{quote}
Note that the continuum hypothesis implies the existence of
a Luzin set, i.e., an uncountable set of reals which meets
every nowhere dense set in a countable set.  Hence, we can think
of Laczkovich's question as asking if one can construct a sort
of weak version of a Luzin set without any extra set theoretic
assumptions.

Recall that a perfect set in a Polish space is a closed nonempty
set without isolated points, and a Polish space is said to be
perfect if it is nonempty and has no isolated points. As we shall
see, the underlying space in the question of Laczkovich can be
taken to be any perfect Polish space. If we ask, as is quite
natural, for the nowhere dense perfect sets in the statement to be
Cantor sets (i.e., sets homeomorphic to the Cantor middle third
set), then we do not know whether the nature of the Polish space
matters.  Even for various standard incarnations of the reals (the
real line, the Baire space, and so on), we have only partial
results on their equivalence in this context. We answered
Laczkovich's question for the Cantor set in 1997 by building a
model where the answer is negative. (And of course the perfect
nowhere dense sets in this case are necessarily Cantor sets.) Very
shortly afterwards, we noticed the more elegant solution suggested
by the first sentence of Remark~\ref{(b) is enough for perf(R) +
model from CS} in which the negative answer is deduced from a
slight variant on a statement proven consistent by Shelah in
\cite{Sh1980}. We show in Section~\ref{Consistency results} that
the stronger conclusion in which, for any perfect Polish space,
the perfect nowhere dense sets can be taken to be Cantor sets
follows from yet another variant on the same statement. The proof
of the consistency of the variants in question is similar to the
proof of Shelah. Unfortunately, the proof is quite technical and
the argument in \cite{Sh1980} is only a brief sketch, so we give
the argument in some detail in Section~\ref{Order-isomorphisms of
everywhere nonmeager sets} in order to be clear. An alternative
model for the negative answer to Laczkovich's question for the
Cantor set is provided by a paper of Ciesielski and Shelah
\cite{CS}. See Remark~\ref{(b) is enough for perf(R) + model from
CS}. In the final section of the paper, we show how a measure
theoretic version of our results can be deduced from results in
Roslanowski and Shelah \cite{RS}. The authors thank Ilijas Farah
for helpful discussions concerning the models constructed in
\cite{RS}.

Write $\perf(X)$ for a Polish space $X$ to mean that for every
nonmeager set $A\sq X$ there is a nowhere dense perfect set $P\sq
X$ such that $A\cap P$ is nonmeager relative to $P$. Write
$\cantor(X)$ if moreover $P$ can be taken to be a Cantor set. Note
that $\perf(X)$ and $\cantor(X)$ are trivially equivalent in
spaces in which nowhere dense perfect sets are necessarily Cantor
sets, e.g., $2^\o$ and $\R$.

We recall for emphasis the following well-known elementary fact of
which we will make frequent use without mention.
\prop{rel nwd in dense subspace}{If $X$ is a topological space and
$Y$ is a dense subspace of $X$, then for any $A\sq Y$, $A$ is
nowhere dense in $Y$ if and only if $A$ is nowhere dense in $X$.
Similarly, $A$ is meager in $Y$ if and only if $A$ is meager in
$X$. \qed}

\section{Relationships between various Polish spaces}
\label{Relationships between various Polish spaces}

We begin by showing that for any two perfect Polish spaces $X$ and
$Y$, $\perf(X)$ and $\perf(Y)$ are equivalent statements.

\prop{equivalence of perf(irr) and perf(Polish)}{
\begin{enumerate}
\item[{\rm(a)}]
Suppose $X$ is a perfect Polish space and $\perf(X)$ holds. Then
$\perf(\o^\o)$ holds.
\item[{\rm(b)}]
$\perf(\o^\o)$ implies $\perf(X)$ for every perfect Polish space
$X$.
\end{enumerate}
 }

\proof  We will use the well-known fact that every perfect Polish
space $X$ has a dense $G_\d$ subset $Y$ homeomorphic to $\o^\o$.
(To get $Y$, first remove the boundaries of the elements of a
countable base for $X$. What remains is a zero-dimensional dense
$G_\d$. Remove a countable dense subset of this dense $G_\d$ and
call the result $Y$. Then $Y$ is a perfect Polish space which is
zero-dimensional and has no compact open sets and hence is
homeomorphic to $\o^\o$.)

(a) Let $Y$ be a residual subspace of $X$ homeomorphic to $\o^\o$.
Let $A$ be a nonmeager set in $Y$. In $X$, $A$ is nonmeager so
there is a nowhere dense perfect set $C$ so that $A$ is nonmeager
in $C$. By replacing $C$ by the closure of one of its nonempty
open subsets, we may assume that $A$ is everywhere nonmeager in
$C$. In particular, $A\cap C$ is dense in $C$. Note that $F =
Y\cap C$ is closed relative to $Y$, is nonempty and has no
isolated points (because it contains $A\cap C$ which is dense in
$C$). Since $F$ is dense in $C$, $A\cap C=F\cap C$ is not meager
in $F$. Also, because $Y$ is dense in $X$ and $F$ is nowhere dense
in $X$, $F$ is also nowhere dense in $Y$.

(b) Let $X$ be a perfect Polish space. Let $Y$ be a residual
subspace of $X$ homeomorphic to $\o^\o$. Let $A\sq X$ be
nonmeager. Then $A\cap Y$ is nonmeager in $X$ and hence in $Y$ as
well since $Y$ is dense. By $\perf(\o^\o)$, there is a nowhere
dense perfect set $C$ in $Y$ such that $A\cap C$ is nonmeager in
$C$.  If $P$ denotes the closure of $C$ in $X$, then, since $C$ is
dense in $P$, $A\cap C$ is nonmeager in $P$ and hence $A\cap P$ is
also nonmeager in $P$. $P$ is perfect since it is the closure of a
nonempty set without isolated points. $P$ is nowhere dense since
it is the closure of a set which is nowhere dense in $Y$ and hence
in $X$ as well.
 \qed

\m

Part (b) holds for $\cantor(\cdot)$ by an easier argument.

\prop{(b) for perf^*}{$\cantor(\o^\o)$ implies $\cantor(X)$ for
every perfect Polish space $X$.}

\proof  Similar to the proof of Proposition~\ref{equivalence of
perf(irr) and perf(Polish)}(b), except that this time the proof
yields a nowhere dense Cantor set $C\sq Y$ such that $A\cap C$ is
nonmeager in $C$ and then we are done.
 \qed

\m

We do not know whether (a) holds for $\cantor(\cdot)$.

\pr{Irrationals vs Cantor set}{Does $\cantor(2^\o)$ imply
$\cantor(\o^\o)$?}

\pr{Unit square vs Cantor set}{Does $\cantor([0,1])$ imply
$\cantor([0,1]\times[0,1])$?}
Of course, $\cantor([0,1])\equiv\perf([0,1]) \equiv\perf(2^\o)
\equiv\cantor(2^\o)$, so these two questions have equivalent
hypotheses.

\m

We introduce one more version of $\perf(X)$ based on the following
observation. If $\perf(\o^\o)$ holds, then for any nonmeager set
$A$, we have a perfect set $P$ such that $A\cap P$ is nonmeager in
$P$. Replacing $P$ by the closure of one of its open sets, we may
assume that $A\cap P$ is everywhere nonmeager in $P$. Then if $P$
has a compact open subset $U$, then $U$ is a Cantor set and $A\cap
U$ is nonmeager in $U$. Otherwise, $P$ itself is homeomorphic to
$\o^\o$. Hence, the perfect set $P$ in the conclusion of
$\perf(\o^\o)$ can always be taken to be either a closed nowhere
dense copy of $\o^\o$ or a Cantor set. Let $\irr(X)$ be the
strengthening of $\perf(X)$ in which we require that the perfect
nowhere dense sets in the definition be homeomorphic to the Baire
space $\o^\o$. Of course a Polish space need not contain any
closed copies of $\o^\o$, so $\irr(X)$ can fail.  However, when
$X=\o^\o$ it would seem reasonable that $\irr(X)$ might hold, and
we will show in the next section that $\irr(\o^\o)$ is indeed
consistent. Its relationship to $\cantor(\o^\o)$ is unclear to us.

\pr{question on irr} {(a) Does $\perf(\o^\o)$ imply that one of
$\irr(\o^\o)$ or $\cantor(\o^\o)$ must hold?  (b) Does either of
$\irr(\o^\o)$ or $\cantor(\o^\o)$ imply the other? }

\section{Consistency results}
\label{Consistency results}

We now turn to the proof of the consistency of $\cantor(\o^\o)$
and $\irr(\o^\o)$. We need a variation on the following result
which forms part of the proof of \cite[Theorem~4.7]{Sh1980} which
states that if ZFC is consistent, then so is ZFC + $\neg$CH +
``There is a universal (linear) order of power $\o_1$.''

\thm{universal order 1}{If ZFC is consistent, then so is ZFC +
both of the following statements.
\begin{enumerate}
\item[{\rm(a)}]
There is a nonmeager set in $\R$ of cardinality $\o_1$.
\item[{\rm(b)}]
Let $A$ and $B$ be everywhere nonmeager subsets of $\R$ of
cardinality $\o_1$. Then $A$ and $B$ are order-isomorphic.
\end{enumerate}
 }

We shall need the following variant of this result.

\thm{universal order 2}{If ZFC is consistent, then so is ZFC +
both of the following statements.
\begin{enumerate}
\item[{\rm(a)$'$}]
Every nonmeager set in $\R$ has a nonmeager subset of cardinality
$\o_1$.
\item[{\rm(b)$'$}]
Let $A$ and $B$ be everywhere nonmeager subsets of $\R$ of
cardinality $\o_1$. Suppose we are given countable dense subsets $A_0\sq A$
and $B_0\sq B$. Then $A$ and $B$ are order-isomorphic by an order
isomorphism taking $A_0$ isomorphically to $B_0$.
\end{enumerate}
 }

\pr{Is added part automatic}{In the presence of (a), does (b)
imply (b)$'$?}

We shall in fact verify in Theorem~\ref{universal order 3} that in
(b)$'$ we can even ask that given pairwise disjoint countable
dense subsets $A_i$, $i<\o$, of $A$ and pairwise disjoint
countable dense subsets $B_i$, $i<\o$, of $B$, the
order-isomorphism of $A$ and $B$ takes $A_i$ isomorphically to
$B_i$ for each $i<\o$. As explained in the introduction, the proof
is similar to the one in \cite{Sh1980}, but as the proof is quite
technical and the argument in \cite{Sh1980} is only a brief
sketch, we need to give the argument in some detail in order to be
clear. We do that in the next section.  Here, we derive the
consequences of interest to us for this paper. The definition of
$\irr(X)$ is given at the end of Section~\ref{Relationships
between various Polish spaces}.

\thm{consequences of interest}{Assume {\rm(a)$'$} and {\rm(b)$'$}.
Then $\cantor(\o^\o)$ and $\irr(\o^\o)$ both hold.}

\proof We will use the following elementary fact.

\faact{extending order isomorphisms}{If $K,L\subseteq\R$ are dense
and $h\colon K\to L$ is an order isomorphism, then $h$ extends to
an order isomorphism of $\R$. \qed}
Suppose that $A\subseteq\R\sm\Q$ is not meager. We wish to find a
Cantor set $C\sq\R\sm\Q$ such that $A\cap C$ is nonmeager relative
to $C$. By (a)$'$, we may assume that $A$ has cardinality exactly
$\o_1$. $A$ is everywhere nonmeager in some open interval
$(a,b)$. Let $C\subseteq (a,b)\sm\Q$ be a Cantor set, and, by
(a)$'$, let $B\subseteq C$ be a set of cardinality $\o_1$ which is
nonmeager relative to $C$. Then $A\cup\Q$ and $A\cup B\cup\Q$ are
both everywhere nonmeager in $(a,b)$ and both have cardinality
$\o_1$. By (b)$'$, there is an order-isomorphism $h\colon
A\cup\Q\to A\cup B\cup\Q$ such that $h[\Q]=\Q$. Extend $h$ to
$(a,b)$ and denote the extension also by $h$. Since $h$ is a
homeomorphism, $h^{-1}[C]$ is a Cantor set and $h^{-1}[B]$ is non
meager relative to $h^{-1}[C]$. Since $h^{-1}[C]\sq\R\sm\Q$ and
$h^{-1}[B]\subseteq A$, we are done.

To get $\irr(\o^\o)$, we make a different choice of $C$ in the
preceding argument.  This time, choose $C$ to be any Cantor set so
that $C\cap\Q$ is dense in $C$. Then $h^{-1}[C]$ will have the
same property, so $h^{-1}[C\sm\Q]=h^{-1}[C]\sm\Q$ is closed
nowhere dense in $\R\sm\Q$ and homeomorphic to $\o^\o$. \qed

\rem{(b) is enough for perf(R) + model from CS}{The reader can
easily verify that a similar but simpler argument yields that
(a)$'$ and (b) imply $\perf(\R)$. An alternative proof of the
consistency of $\perf(\R)$ can by had by using Theorem 2 of
\cite{CS} which states that the following statement is consistent
relative to ZFC:
\begin{quote}
For every $A\subseteq 2^\omega\times 2^\omega$ for which the sets
$A$ and $A^c=(2^\omega\times 2^\omega\setminus A)$ are nowhere
meager in $2^\omega\times 2^\omega$ there is a homeomorphism
$f:2^\omega\to 2^\omega$ such that the set $\{x\in 2^\omega:
(x,f(x))\in A\}$ does not have the Baire property in $2^\omega$.
\end{quote}
Note that the map $2^\o\to f$ given by $x\mapsto(x,f(x))$ is a
homeomorphism. Hence the conclusion could be stated as ``$f\cap A$
does not have the Baire property in $f$''. Since $2^\omega\times
2^\omega$ is homeomorphic to $2^\omega$ and the graph of a
homeomorphism of $2^\o$ is a perfect nowhere dense set in
$2^\omega\times 2^\omega$, the statement above implies the
following special case of $\perf(2^\o)$ (which is equivalent to
$\perf(\R)$).
\begin{quote}
For every $A\subseteq 2^\omega$ for which the sets $A$ and $A^c$
are both nowhere meager in $2^\omega$, there is a perfect nowhere
dense set $P$ such that the set $A\cap P$ is not meager in $P$.
\end{quote}
To reduce $\perf(2^\o)$ to this special case, consider a nonmeager
set $A\subseteq 2^\omega$.  $A$ is everywhere nonmeager in some
clopen set $U$. If $A$ is comeager in some clopen set, then it
contains a nowhere dense perfect set and we are done. Hence we may
assume that, relative to $U$, $A$ and $A^c$ are both everywhere
nonmeager. The clopen set $U$ is homeomorphic to $2^\omega$, so we
now find ourselves in the special case described above.
 }

\section{Order-isomorphisms of everywhere nonmeager sets}
\label{Order-isomorphisms of everywhere nonmeager sets}

We now turn to the proof of the consistency of (a)$'$
and (b)$'$. We begin by recalling the basic properties
of oracle-cc forcing. See \cite[Chapter IV]{Sh1998}
for the details. A version of this material is also
explained in \cite[Sections 4--6]{Bu}.

\defi{oracle}{A sequence
\[
\barM=\la M_\d:\d\ \mbox{is a limit ordinal}\ <\o_1\ra
\]
is called an {\em oracle} if each $M_\d$ is a countable transitive
model of a sufficiently large fragment of ZFC, $\d\in M_\d$ and
for each $A\sq\o_1$, $\{\d:A\cap\d\in M_\d\}$ is stationary in
$\o_1$.}
The meaning of ``sufficiently large'' depends on the context. In a
particular proof, some fragment of ZFC for which models can be
produced in ZFC must suffice for all the oracles in the proof. The
existence of an oracle is equivalent to $\diamondsuit$, (see
\cite[Theorem II 7.14]{Ku}) and hence implies CH. We limit the
definition of the $\barM$-chain condition to partial orders of
cardinality $\o_1$. This covers our present needs.

\m

Associated with an oracle $\barM$, there is a filter $\Trap$
generated by the sets
\[
\{\d<\o_1:\d\ \mbox{is a limit ordinal and}\ A\cap\d\in
M_\d\},\quad A\sq\o_1.
\]
This is a proper normal filter containing all closed unbounded
sets.

\defi{completely embedded modulo D}{
If $P$ is any partial order, $P'\sq P$, and ${\mathfrak D}$ is any
class of sets, then we write $P'<_{\mathfrak D}P$ to mean that
every predense subset of $P'$ which belongs to ${\mathfrak D}$ is
predense in $P$.}

\defi{oracle-cc}{A partial order $P$ {\it satisfies the $\barM$-chain
condition}, or simply {\it is $\barM$-cc}, if there is a
one-to-one function $f\colon P\to\o_1$ such that
\[
\{\d<\o_1:\d\ \mbox{is a limit ordinal and}\
f^{-1}(\d)<_{M_{\d,f}}P\}
\]
belongs to $\Trap$, where $M_{\d,f}=\{f^{-1}(A):A\sq\d,\,A\in
M_\d\}$.}

It is not hard to verify that if $P$ is $\barM$-cc, then $P$ is
ccc. Also, any one-to-one function $g\colon P\to\o_1$ can replace
$f$ in the definition.

\prop{oracle-cc properties}{The $\barM$-cc satisfies the following
properties.
\begin{enumerate}
\item[{\rm(1)}]
If $\a<\o_2$ is a limit ordinal, $\la\la
P_\beta\ra_{\beta\leq\a},\la \dot Q_\beta\ra_{\beta<\a}\ra$ is a
finite-support $\a$-stage iteration of partial orders, and for
each $\beta<\a$, $P_\beta$ is $\barM$-cc, then $P_\a$ is
$\barM$-cc.

\item[{\rm(2)}]
If $P$ is $\barM$-cc, then there is a $P$-name $\barM^*$ for an
oracle such that for each $P$-name $\dot Q$ for a partial order,
if $\Vdash_P$~``$\dot Q$ is $\barM^*$-cc'' then $P*\dot Q$ is
$\barM$-cc.

\item[{\rm(3)}]
If $\barM_\a$, $\a<\o_1$, are oracles, then there is an oracle
$\barM$ such that for any partial order $P$, if $P$ is $\barM$-cc,
then $P$ is $\barM_\a$-cc for all $\a<\o_1$.
\end{enumerate}
 }

We will need the following lemmas.

\lem{Main Lemma 1}{Let $\barM=\la M_\d:\d<\o_1\ra$ be an oracle
and let $A$ and $B$ be everywhere nonmeager subsets of $\R$.
Suppose we are given pairwise disjoint countable dense subsets
$A_i$, $i<\o$, of $A$ and pairwise disjoint countable dense
subsets $B_i$, $i<\o$, of $B$. Then there is a forcing notion $P$
satisfying the $\barM$-cc such that for every $G\sq P$ generic
over $V$, $V[G]\models$  $A$ and $B$ are order-isomorphic by an
order isomorphism taking $A_i$ isomorphically to $B_i$ for each
$i<\o$.
 }

\proof Fix well-orderings of $A$ and $B$ in type $\o_1$. (CH holds
because there is an oracle.) We will inductively define one-to-one
enumerations $\la a_\a:\a<\o_1\ra$ of $A$ and $\la
b_\a:\a<\o_1\ra$ of $B$ and functions $f_\d$, $\d<\o_1$. We let
$A_\d=\{a_\a:\o\d\leq\a<\o(\d+1)\}$ and
$B_\d=\{b_\a:\o\d\leq\a<\o(\d+1)\}$ for $\d<\o_1$. For $A'\sq A$
and $B'\sq B$, let $P(A',B')$ denote the set of finite partial
order-preserving maps $p\colon A'\to B'$ such that $p[A_\d]\sq
B_\d$ for all $\d<\o_1$. We also use the notation
\[
A\restr\a=\{a_\beta:\beta<\a\},\ B\restr\a=\{b_\beta:\beta<\a\}.
\]
We will arrange that the following conditions hold.
\begin{enumerate}
\item[(1)] The sets $A_\d$ and $B_\d$ are dense in $\R$.
\item[(2)]
For $\d<\o$, the sets $A_\d$ and $B_\d$ are as in the hypothesis.
\item[(3)]
For each $\d<\o_1$, $f_\d$ is a bijective map of
$P(A\restr\o\d,B\restr\o\d)$ onto $\o\d$.
\item[(4)] For each $\d<\d'<\o_1$, $f_\d\sq f_{\d'}$.
\item[(5)]
For each infinite $\d<\o_1$, the predense subsets of
$P(A\restr\o\d,B\restr\o\d)$ which have the form $f_\d^{-1}[S]$
for some $S\in \bigcup_{\eta\leq\d}M_{\eta}$ remain predense in
$P(A\restr\o(\d+1),B\restr\o(\d+1))$.
\end{enumerate}
To do this, we proceed as follows. The construction of the
functions $f_\d$ is dictated by (4) at limit stages, and
$f_{\d+1}$ is an arbitrary extension of $f_\d$ satisfying (3). The
elements of $A_\d$ and $B_\d$ for $\d<\o$ are given by (2). For
$\d\geq\o$, by induction on $\d$ we choose the elements of $A_\d$
and $B_\d$ by alternately defining $a_{\o\d+n}$ and $b_{\o\d+n}$,
beginning with $a_{\o\d}$ when $\d$ is even and with $b_{\o\d}$
when $\d$ is odd. Let us illustrate the construction with the case
where $\d$ is even. Fix an enumeration $\la I_m:0<m<\o\ra$ of the
nonempty open intervals with rational endpoints. The first element
$a_{\o\d}$ is simply the least element, in the well-ordering of
$A$ fixed at the beginning of the proof, which is different from
any of the elements of $A$ chosen so far. We now choose
$b_{\o\d},a_{\o\d+1}, b_{\o\d+1},a_{\o\d+2},b_{\o\d+2},\dots$ in
that order. For $n>0$, we pick $a_{\o\d+n}$ and $b_{\o\d+n}$ from
$I_n$ to ensure $A_\d$ and $B_\d$ will be dense.

To choose one of these elements, say $b_{\o\d+n}$, let $N$ be a
countable elementary submodel of $H_\t$, for a suitably large
$\t$, such that $A$, $B$, $f_\d$, the sequences $\la
a_\a:\a\leq\o\d+n\ra$ and $\la b_\a:\a<\o\d+n\ra$, and
$\bigcup_{\eta\leq\d}M_\eta$ are all elements of $N$. Choose
$b_{\o\d+n}$ to be a member of $B$ which is a Cohen real over $N$.

We must check that the construction gives (5). Let $D$ be a
predense subset of $P(A\restr\o\d,B\restr\o\d)$ of the appropriate
form. In particular, we have $D\in N$. We will show that $D$
remains predense in $P(A\restr\o\d+n+1,B\restr\o\d+n+1)$.

\rem{why it works}{We are showing by induction on $n$ that $D$
remains predense in $P(A\restr\o\d+n+1,B\restr\o\d+n)$ and then in
$P(A\restr\o\d+n+1,B\restr\o\d+n+1)$. (This establishes (5) since
each member of $P(A\restr\o(\d+1),B\restr\o(\d+1))$ belongs to
$P(A\restr\o\d+n,B\restr\o\d+n)$ for some $n<\o$.) Our current
stage has the second form. Note that at the stage $n=0$, we first
consider the passage from $P(A\restr\o\d,B\restr\o\d)$ to
$P(A\restr\o\d+1,B\restr\o\d)$. But these two partial orders are
equal because there is no legal target value for $a_{\o\d}$ until
$b_{\o\d}$ is chosen. So the preservation of the predense sets
trivially holds at that stage. In particular, it does not matter
that $a_{\d\o}$ is not Cohen generic over the previous
construction.}
Let
\[
p\in P(A\restr\o\d+n+1,B\restr\o\d+n+1)\sm
P(A\restr\o\d+n+1,B\restr\o\d+n).
\]
Then $p$ has the form $q\cup\{(a,b_{\o\d+n})\}$ for some $q\in
P(A\restr\o\d+n+1,B\restr\o\d+n)$ and $a\in\{a_{\o\d+m}:m\leq
n\}$. Fix $r\in D$. The set
\[
\{b\in\R: q\cup\{(a,b)\}\ \mbox{is compatible with}\ r\}\in N
\]
is open and hence its complement $C_r$ is closed, as is the set
$C_D=\bigcap_{r\in D} C_r$ of $b$ for which $q\cup\{(a,b)\}$ is
incompatible with every member of $D$. Since $p$ is an partial
order isomorphism, there are open rational intervals $J_1$ and
$J_2$ such that $J_1\cap\dom p=\{a\}$, $J_2\cap \range
p=\{b_{\o\d+n}\}$. Note that whenever $x\in J_1$ and $b\in J_2$,
$q\cup\{(x,b)\}$ is a partial isomorphism.

\clm{C_D is nowhere dense in J_2}{$C_D$ is nowhere dense in
$J_2$.}

\proof Fix a nonempty open subinterval $J$ of $J_2$. There is an
extension of $q$ by members of $A_0\times B_0$---the point of
using $A_0$ and $B_0$ being simply that they are dense and
contained in $A(\o\d+n+1)$ and $B(\o\d+n)$, respectively---which
adds two points in $J_1\times J$ straddling the line $x=a$. So
this extension has the form
\[
q'=q\cup\{(x_1,y_1),(x_2,y_2)\},\ x_1<a<x_2,\,y_1<y_2
\]
where $(x_1,x_2)\sq J_1$ and $(y_1,y_2)\sq J$. Since $q'\in
P(A\restr\o\d+n+1,B\restr\o\d+n)$, by the induction hypothesis $D$
must have an element $r$ compatible with this extension. Since
$a\not\in A(\o\d)$, we have $a\not\in\dom r$. Let $x'_1$, $x'_2$
be the closest members of $\dom(q'\cup r)$ to the left and right
of $a$, respectively. Write $y'_1=r(x'_1)$, $y'_2=r(x'_2)$. Then
$(y'_1,y'_2)\sq (y_1,y_2)\sq J$ and for any choice of $b\in
(y'_1,y'_2)$, $q\cup\{(a,b)\}$ is compatible with $r$. Hence
$(y'_1,y'_2)$ is disjoint from $C_r$ and hence from $C_D$. This
proves the claim.
 \qed

\m

\n Thus, $b_{\o\d+n}\not\in C_D$ and hence $p$ is compatible with
some member of $D$. This establishes (5).

\m

Now take $P=P(A,B)$. The fact that $P$ forces the desired
order-isomorphism of $A$ and $B$ is clear from (1) and (2). To see
that $P$ is $\barM$-cc, let $f=\bigcup_{\d<\o_1}f_\d\colon
P\to\o_1$. For any $\d<\o_1$ we have
$f^{-1}[\o\d]=P(A(\o\d),B(\o\d))$ and for each $S\sq\o\d$ whenever
a set $D$ of the form $f^{-1}[S]=f_\d^{-1}[S]$ belongs to $M_\d$
and is predense in $P(A(\o\d),B(\o\d))$, a simple induction on
$\d'$ using (5) shows that if $\d$ is infinite and
$\d<\d'\leq\o_1$, then $D$ is predense in $P(A(\o\d'),B(\o\d'))$.
In particular, $D$ is predense in $P=P(A(\o_1),B(\o_1))$. For a
club of $\d<\o_1$ we have $\o\d=\d$, so this shows that $P$
satisfies the $\barM$-cc.
 \qed

\lem{Main Lemma 2}{Assume $\diamondsuit$. Let $A$ be a nonmeager
subset of $\R$. Then there is an oracle $\barM=\la
M_\d:\d<\o_1\ra$ such that if $P$ is any partial order satisfying
the $\barM$-cc, then $\Vdash_P$ ``$A$ is nonmeager''.
 }

\proof This is \cite[Example IV 2.2]{Sh1998}.
 \qed

\thm{universal order 3}{If ZFC is consistent, then so is ZFC +
both of the following statements.
\begin{enumerate}
\item[{\rm(a)}]
Every nonmeager set in $\R$ has a nonmeager subset of cardinality
$\o_1$.
\item[{\rm(b)}]
Let $A$ and $B$ be everywhere nonmeager subsets of $\R$ of
cardinality $\o_1$. Suppose we are given pairwise disjoint
countable dense subsets $A_i$, $i<\o$, of $A$ and pairwise
disjoint countable dense subsets $B_i$, $i<\o$, of $B$. Then $A$
and $B$ are order-isomorphic by an order isomorphism taking $A_i$
isomorphically to $B_i$ for each $i<\o$.
\end{enumerate}
 }

\proof Start with a ground model of $V=L$. Fix a diamond sequence
\[
\la (f_\a,g_\a,h_\a):\a<\o_2,\,\cof(\a)=\o_1\ra
\]
for trapping triples $(f,g,h)$ consisting of:
\begin{enumerate}
\item[(1)]
A function $f\colon\o_2\to([\o_2]^{\leq\o})^\o$. The idea of $f$
is that, with $\o_2$ identified with the ccc partial order we are
about to build, $[\o_2]^{\leq\o}$ contains the maximal antichains.
Thus, $([\o_2]^{\leq\o})^\o$ contains a name for each real number
(construed as a subset of $\o$). Then for any nonmeager set $X$ in
the extension, we can find a ground model function
$f\colon\o_2\to([\o_2]^{\leq\o})^\o$ enumerating the names of the
elements of $X$.
\item[(2)]
Functions $g,h\colon\o_1\to([\o_2]^{\leq\o})^\o$ intended to
represent (enumerations of the names for the elements of)
everywhere nonmeager sets of cardinality $\o_1$ with each of the
sets $\{g(\o i+n):n<\o\}$ and $\{h(\o i+n):n<\o\}$, for $i<\o$,
dense in $\R$.
\end{enumerate}
So for each $\a<\o_2$ of cofinality $\o_1$,
$f_\a\colon\a\to([\a]^{\leq\o})^\o$, and
$g_\a,h_\a\colon\o_1\to([\a]^{\leq\o})^\o$. Also, for each
$(f,g,h)$ as in (1)--(2), $\{\a<\o_2:\cof(\a)=\o_1$,
$f\restr\a=f_\a$, $g\restr\a=g_\a$ and $h\restr\a=h_\a\}$ is
stationary in $\o_2$.

We will inductively define an $\o_2$-stage finite support
iteration
\[
\la\la P_\a\ra_{\a\leq\o_2},\la\dot Q_\a\ra_{\a<\o_1}\ra
\]
as well as a $P_\a$-names $\barM_\a$ for oracles and one-to-one
functions $F_\a\colon P_\a\to\o_2$ for $\a<\o_2$ such that the
range of each $F_\a$ is an initial segment of $\o_2$ which
includes $\a$ and for $\beta<\a<\o_2$, we have $F_\beta\sq F_\a$.
(At each stage, $F_\a$ is any function satisfying these
conditions.)

For $\a<\o_2$, we will let $\dot X_\a$ denote the $P_\a$-name for
the set of real numbers whose elements have the names
\[
{\textstyle\bigcup_{n<\o}}\{n\}\times
F^{-1}_\a(f_\a(\xi)(n)),\quad \xi<\a.
\]
Similarly, we will let $\dot A_\a$ and $\dot B_\a$ denote the
$\o_1$-sequences of $P_\a$-names for real numbers
\[
\left\langle{\textstyle\bigcup_{n<\o}}\{n\}\times
F^{-1}_\a(g_\a(\xi)(n)):\xi<\o_1\right\rangle
\]
and
\[
\left\langle{\textstyle\bigcup_{n<\o}}\{n\}\times
F^{-1}_\a(h_\a(\xi)(n)):\xi<\o_1\right\rangle
\]
respectively. At stage $\a<\o_2$ of the construction, if
$\cof(\a)=\o_1$ and if
\[
\Vdash_{P_\a}\dot X_\a\ \mbox{is not meager},
\]
then we use Lemma~\ref{Main Lemma 2} to get a $P_\a$-name
$\barM'_\a$ for an oracle so that if $P$ is any forcing notion
which satisfies the $\barM'_\a$-cc, then $X_\a$ remains nonmeager
after forcing with $P$. Otherwise, in particular if
$\cof(\a)\not=\o_1$, we let $\barM'_\a$ be any $P_\a$-name for an
oracle.

For $\beta<\a$, let $P_{\beta\a}$ be the usual $P_\beta$-name for
a partial order such that $P_\a$ is isomorphic to a dense subset
of $P_\beta*P_{\beta\a}$ (see [Ba]). Let $\barM_{\beta\a}$ be a
$P_\a$-name for an oracle such that
$$
\mbox{If}\ \Vdash_{P_{\beta}}\ \mbox{``}P_{\beta,\a}\ \mbox{is}\
\barM_\beta\mbox{-cc and}\ \Vdash_{P_{\beta,\a}}\dot Q_\a\
\mbox{is}\ \barM_{\beta\a}\mbox{-cc''},
$$

\vspace{-10pt} 

\n(1)

\vspace{-15pt}

$$
\mbox{then}\ \Vdash_{P_{\beta}}\
\mbox{``}P_{\beta,\a+1}=P_{\beta,\a}*\dot Q_\a\ \mbox{is}\
\barM_\beta\mbox{-cc}\mbox{''}.
$$
(There is such an $\barM_{\beta\a}$ by Proposition~\ref{oracle-cc
properties}(2). In (1), $\barM_{\beta\a}$ is actually a
$P_{\beta}$-name for a $P_{\beta,\a}$-name for an oracle. We
denote the corresponding $P_\a$-name also by $\barM_{\beta\a}$.)

Let $\barM_\a$ be a $P_\a$-name for an oracle such that
$$
\Vdash_{P_\a}\ \mbox{``If}\ \dot Q_\a\ \mbox{is}\
\barM_\a\mbox{-cc, then} \ \dot Q_\a\ \mbox{is}\
\barM'_\a\mbox{-cc and}\ \barM_{\beta\a}\mbox{-cc for all}\
\beta<\a\mbox{''.}\leqno{(2)}
$$
(Use Proposition~\ref{oracle-cc properties}(3).)

Now, if $\cof(\a)=\o_1$ and if
\[
\Vdash_{P_\a}\ \mbox{The ranges of}\ \dot A_\a,\dot B_\a\
\mbox{are everywhere nonmeager and each of the sets}
\]
\[
\hspace{45pt}\{\dot A_\a(\o i+n):n<\o\},\,\{\dot B_\a(\o
i+n):n<\o\},\ \mbox{for}\ i<\o,\ \mbox{is dense in}\ \R.
\]
then use Lemma~\ref{Main Lemma 1} to get a $P_\a$-name $\dot Q_\a$
for a partial order satisfying the $\barM_\a$-cc and forcing an
isomorphism between $A_\a$ and $B_\a$ as described in the
statement of the lemma. In all other cases, take $\dot Q_\a$ to
name the partial order $Q$ for adding one Cohen real. We have thus
$$
\Vdash_{P_\a}\ \mbox{``}\dot Q_\a\ \mbox{satisfies the
$\barM_\a$-cc''}.\leqno{(3)}
$$
Now suppose that for some $P_{\o_2}$-name $\dot X$ we have
\[
\Vdash_{P_{\o_2}}\dot X\ \mbox{is not meager}.
\]
(Every nonmeager set in any extension has a name forced by the
weakest condition to be nonmeager since there is always a
nonmeager set.) Fix a name $\dot f$ such that
\[
\Vdash_{P_{\o_2}}\dot f\colon\o_2\to\dot X\ \mbox{is onto}.
\]
Then define $f\colon\o_2\to([\o_2]^{\leq\o})^\o$ so that if
\[
\tau_\xi={\textstyle\bigcup_{n<\o}}\{n\}\times F^{-1}(f(\xi)(n)),\
\xi<\o_2,
\]
then for each $\xi<\o_2$,
\[
\Vdash_{P_{\o_2}}\dot f(\xi)=\tau_\xi.
\]
There is a closed unbounded set $C\sq\o_2$ such that for each
$\a\in C$ of cofinality $\o_1$ we have:
\begin{enumerate}
\item[(i)] $f\restr\a\colon\a\to([\a]^{\leq\o})^\o$.
\item[(ii)] $\forall\xi<\a$, $\tau_\xi$ is a $P_\a$-name.
\item[(iii)] $\Vdash_{P_\a}\{\tau_\xi:\xi<\a\}$ is not meager.
\end{enumerate}
(For (iii), note that when $\a$ has cofinality $\o_1$, each
$P_\a$-name for a meager set is a $P_\beta$-name for some
$\beta<\a$. Thus, if $M$ is an elementary submodel of $H_\t$ for a
suitably large $\t$ such that $|M|=\o_1$, $M^\o\sq M$,
$\la\tau_\xi:\xi<\o_2\ra\in M$ and $\a=M\cap\o_2\in\o_2$ has
cofinality $\o_1$, then for each (nice) $P_\a$-name $\s$ for a
meager Borel set set, we have $\s\in M$ and hence $M$ knows about
a maximal antichain of conditions each deciding a $\xi$ for which
$\tau_\xi$ is forced not to be in $\s$. The antichain is countable
and hence contained in $M$. For each condition in the antichain,
the least $\xi$ which it decides is in $M$ and hence below $\a$.
Hence $\Vdash_{P_{\a}}$ ``$\{\tau_\xi:\xi<\a\}$ is not contained
in $\s$''.)

Choose such an $\a$ of cofinality $\o_1$ for which
$f\restr\a=f_\a$. By (i) and (ii), the definition of $\tau_\xi$
would not change if we used $f_\a$ instead of $f$ and $F_\a$
instead of $F$. Then from the definition of $\dot X_\a$ we get
\[
\Vdash_{P_{\a}}\dot X_\a=\{\tau_\xi:\xi<\a\}.
\]
So at stage $\a$ we chose a $P_\a$-name $\barM_\a$ and we arranged
that
\[
\Vdash_{P_\a}\ \mbox{``}P_{\a,\g}\ \mbox{is}\
\barM_\a\mbox{-cc''.}
\]
(This follows easily by induction on $\g\geq\a$ and
Propositions~\ref{oracle-cc properties}(1,2). (Recall that
$P_{\a,\g}$ can be viewed in canonical way as an iteration: see
\cite{Ba}. At limits $\g$ use Propositions~\ref{oracle-cc
properties}(1). At stages $\g+1$, use (3) to get $\Vdash_{P_\g}$
``$\dot Q_\g$ satisfies the $\barM_\g$-cc'' and then use (2) and
(1) with $(\beta,\a)$ replaced by $(\a,\g)$.)

Hence, by the choice of $\barM_\a$,
$$
\Vdash_{P_\a}\Vdash_{P_{\a,\g}}\dot X_\a\ \mbox{is not
meager}\leqno{(4)}
$$
from which it follows that
$$
\Vdash_{P_\a}\Vdash_{P_{\a,\o_2}}\dot X_\a\ \mbox{is not meager}
$$
since if this failed then we would have
\[
p\Vdash_{P_\a} q\Vdash_{P_{\a,\o_2}}\dot X_\a\sq \dot B
\]
for some conditions $p\in P_\a$, $q\in P_{\a,\o_2}$ and some name
$\dot B$ for a meager Borel set. But then for some $\g$, we have
$\a<\g<\o_2$, $q\in P_{\a,\g}$ and $\dot B$ is a $P_\g$-name and
this contradicts (4).

By what we have established, there are guaranteed to be sets of
cardinality $\o_1$ which are not meager in any extension by
$P_{\o_2}$. Hence there are guaranteed to be everywhere nonmeager
sets of cardinality $\o_1$. Suppose that for some $P_{\o_2}$-names
$\dot A$ and $\dot B$ for $\o_1$-sequences we have
\[
\Vdash_{P_{\o_2}}\ \mbox{The ranges of}\ \dot A,\dot B\ \mbox{are
everywhere nonmeager and each of the sets}
\]
\[
\hspace{45pt}\{\dot A(\o i+n):n<\o\},\,\{\dot B(\o i+n):n<\o\},\
\mbox{for}\ i<\o,\ \mbox{is dense in}\ \R.
\]
(By what we just said, every pair of everywhere nonmeager sets $A$
and $B$ of cardinality $\o_1$, together with choices of countably
many disjoint countable dense subsets of each one, has a name such
that the weakest condition forces the desired properties.)

Define $g,h\colon\o_1\to([\o_2]^{\leq\o})^\o$ so that if
\[
\s_\xi={\textstyle\bigcup_{n<\o}}\{n\}\times
F^{-1}(g(\xi)(n)),\quad \xi<\o_1
\]
and
\[
\tau_\xi={\textstyle\bigcup_{n<\o}}\{n\}\times
F^{-1}(h(\xi)(n)),\quad \xi<\o_1
\]
then for each $\xi<\o_1$,
\[
\Vdash_{P_{\o_2}}\dot A(\xi)=\s_\xi
\]
and
\[
\Vdash_{P_{\o_2}}\dot B(\xi)=\tau_\xi.
\]
For all large enough $\a<\o_2$, we have:
\begin{enumerate}
\item[(i)]
$g,h\colon\o_1\to([\a]^{\leq\o})^\o$.
\item[(ii)]
$\forall\xi<\a$, $\s_\xi$ and $\tau_\xi$ are $P_\a$-names.
\end{enumerate}
Choose any such $\a$ of cofinality $\o_1$. By (i) and (ii), the
definitions of $\s_\xi$ and $\tau_\xi$ would not change if we used
$g_\a$ instead of $g$,  $h_\a$ instead of $h$,  and $F_\a$ instead
of $F$. Then from the definitions of $\dot A_\a$ and $\dot B_\a$
we get
\[
\Vdash_{P_\a}\ \mbox{The ranges of}\ \dot A_\a,\dot B_\a\
\mbox{are everywhere nonmeager and each of the sets}
\]
\[
\hspace{45pt}\{\dot A_\a(\o i+n):n<\o\},\,\{\dot B_\a(\o
i+n):n<\o\},\ \mbox{for}\ i<\o,\ \mbox{is dense in}\ \R.
\]
(Being everywhere nonmeager is trivially downward absolute.) Then
$\dot Q_\a$ was chosen to add an order isomorphism between $A_\a$
and $B_\a$ of the desired type.

This completes the proof of the theorem.
 \qed

\section{A measure-theoretic analog of $\perf(2^\o)$}
\label{measure-theoretic analog}

A measure theoretic version of Laczkovich's question is not
completely obvious because perfect sets carry many measures.  We
consider the following measures on $2^\o$ which we will call
canonical. Given $P\subseteq 2^\omega$ a perfect set, define
$$T_P=\{s\in 2^{<\omega}: P\cap [s]\not=\emptyset\}$$
We say that $s\in T_P$ splits iff both $s0$ and $s1$ are in $T_P$.
The canonical measure $\mu_P$ is the one supported by $P$ and
determined by declaring $\mu_P([s])=1/{2^n}$ iff $s\in T_P$ and
$|\{i<|s|:s\res i$ splits$\}|=n$. An equivalent view is to take
the natural map from $2^{<\omega}$ to the splitting nodes of $T_P$
and the homeomorphism $h:2^\omega\to P\sq 2^\o$ induced by it and
then $\mu_P$ is the measure corresponding to the product measure
$\mu$ on $2^\omega$, i.e., $\mu_P(A)=\mu(h^{-1}(A))$.

\thm{measure analog}{It is relatively consistent with ZFC that for
any set $B\subseteq 2^\omega$ which is not of measure zero, there
exists a perfect set $P$ of measure zero such that $B\cap P$ does
not have measure zero in the canonical measure $\mu_P$ on $P$.}

\proof The model is the one used by Ros\l anowski and Shelah in
the proof of \cite[Theorem~3.2]{RS}. It is obtained by forcing
over a model of CH with an $\o_2$-stage countable support
iteration $\la\la\poset_\a\ra_{\a\leq\o_2},\la\dot{\mathbb
Q}_\a\ra_{\a<\o_2}\ra$ of the measured creature forcing ${\mathbb
Q}={{\mathbb Q}}^{\rm mt}_4(K^*, \Sigma^*,{\bf F}^*)$ defined in
\cite[Section~2]{RS}. We use the notation of \cite{RS} concerning
this partial order. The definition involves in particular a
rapidly growing sequence of powers of $2$, $\la
N_i=2^{M_i}:i<\o\ra$. Forcing with ${\mathbb Q}$ gives rise to a
continuous function $h\colon \prod_{i<\o}N_i\to 2^\o$.  We will
make use of the following result concerning this function. The
measure on $\prod_{i<\omega}N_i$ in this proposition is the
product of the uniform probability measures on the factors and the
measure on $2^\o$ is the usual product measure. In the remainder
of this proof, we denote both of these measures, as well as their
product, by $\mu$, letting the context distinguish them.
\begin{quote}\cite[Proposition~2.6]{RS}
Suppose that $A\subseteq \prod_{i<\omega}N_i \times 2^\omega$ is a
set of outer measure one. Then, in $V^{\mathbb Q}$, the set
$$\{x\in\prod_{i<\omega}N_i : (x,h(x))\in A\}$$
has outer measure one.
\end{quote}
We shall also need to know that ${\mathbb Q}$ is proper and that
countable support iterations of ${\mathbb Q}$ preserve Lebesgue
outer measure. The former is \cite[Corollary~1.14]{RS}. The latter
is explained in the proof of \cite[Theorem~3.2]{RS}. (The
explanation refers the reader to some very general preservation
theorems for iterated forcing. For the reader who wants to verify
this without learning these general theorems, we indicate that it
also follows from the special case of these theorems, preservation
of $\sqsubset^{\rm random}$, given in \cite{Go 1993} by imitating
the proof in \cite{Pa 1996} that Laver forcing satisfies what is
called there $\bigstar$ and by noting that $\bigstar$ implies
preservation of $\sqsubset^{\rm random}$.)

\smallskip

Recall that $N_i=2^{M_i}$. We identify $N_i$ with the set of
binary sequences of length $M_i$. The map
$h:\prod_{i<\omega}N_i\to 2^\omega$ is determined from a generic
sequence of finite maps $(W(i):N_i\to 2: i<\omega)$ added by
${\mathbb Q}$. $h$ is defined by $h(x)(i)=W(i)(x(i))$ for each
$i$. We use the $W(i)$'s to define a perfect set $P\subseteq
2^\omega$ by the condition that $x\in P$ if and only if there
exists $y \in 2^\omega$ such that $x$ is the concatenation of the
sequence $s_0,i_0,s_1,i_1,\ldots$ where $y$ is the concatenation
of $s_0,s_1,s_2,\ldots$  and where each $s_k$ has length $M_k$ and
$i_k=W(k)(s_k)\in \{0,1\}$. $P$ is essentially the same as the
graph of $h$ but we spell out the details to be sure the canonical
measure is the one we want.  Another way to define $P$ is as
follows:
\begin{enumerate}
\item[(i)]
Let $l_i=M_i+\sum_{j<i} (M_j+1)$.  Let $l_{-1}=-1$. The $l_i$,
$i<\o$, are the nonsplitting levels of the tree $T_P$ which can be
determined by the next two conditions.

\item[(ii)]
If $s\in T_P$ and $l_{i-1}<|s|<l_i$ then both $s0$ and $s1$ are in
$T_P$.

\item[(iii)]
If $s\in T_P$ and $|s|=l_i$, then only $sj$ in $T_P$ where
$W(i)(t)=j$ and $s=rt$ is the concatenation of $r$ and $t$ where
$|t|=M_i$ and $r$ has the appropriate length.

\item[(iv)]
Define
$$
P=[T_P]\ \stackrel{\mbox{\scriptsize def}}{=}\ \{x\in 2^\omega:
\forall n\;\; x\res n\in T_P\}
$$
Every time we pass a nonsplitting level $l_i$ we lose half the
measure and so $P$ is a perfect set of measure zero for the usual
measure on $2^\o$.
\end{enumerate}
Let $\rho:\prod_{i<\omega}N_i \times 2^\omega\to 2^\omega$ be the
natural homeomorphism given by  $\rho(x,z)$ is the concatenation
of the sequence $x_0,z_0,x_1,z_1,\ldots$ where we are identifying
$N_i$ with the set of binary sequences of length ${M_i}$.

\clm{rho is measure preserving}{$\rho$ is measure-preserving.}

\proof By a standard uniqueness theorem for the extension of a
measure from an algebra to the $\s$-algebra it generates, it
suffices to verify that $\rho^{-1}[C]$ has the same measure as $C$
for every clopen set $C\sq 2^\o$.  Every clopen set $C\sq 2^\o$
can be partitioned into clopen sets of the form $[r]$, where for
some $k<\o$, $s^r\in\prod_{i<k}N_i$, $t^r\in 2^k$ and $r$ is the
concatenation of $s^r_0,t^r_0,\dots,s^r_{k-1}t^r_{k-1}$. (These
are simply the basic open sets $[r]$ for which $r$ has length
$\sum_{i<k}(M_i+1)$ for some $k<\o$.) Hence it suffices to verify
$\mu(\rho^{-1}[[r]])=\mu([r])$ for $r$ of this form. We have
\[
\textstyle{\mu(\rho^{-1}[[r]])=\mu([s^r]\times[t^r])=\left(
\prod_{i<k}2^{-M_i}\right)2^{-k}=
2^{-\sum_{i<k}(M_i+1)}=\mu([r]).}
\]
This proves the claim.
 \qed

\m

\n Let $g\colon \prod_{i<\omega}N_i\to \prod_{i<\omega}N_i\times
2^\o$ be the homeomorphism of $\prod_{i<\omega}N_i$ onto the graph
of $h$ given by $g(x)=(x,h(x))$. We have $\rho[h]=P$ (i.e., the
graph of $h$ corresponds to $P$ under $\rho$).

\clm{identification of measures}{For any Borel set $B\subseteq
2^\omega$
$$\mu_P(B)=\mu \left(g^{-1}[\rho^{-1}[B]]\right)$$
and similarly for outer measure.}

\proof Since the range of $g$ is the graph of $h$, we have
\[
\mu\left( g^{-1}[\rho^{-1}[B]]\right)=\mu\left(
g^{-1}[\rho^{-1}[B]\cap h]\right)=\mu\left( g^{-1}[\rho^{-1}[B\cap
P]]\right).
\]
Similarly, since $\mu_P$ concentrates on $P$,
$\mu_P(B)=\mu_P(B\cap P)$. Hence, it suffices to prove the claim
for Borel subsets of $P$.

Given $s\in\prod_{i<k}N_i$, define $t^s\in 2^k$ by
$t^s_i=W(i)(s(i))$ for all $i<k$, and write
$r^s=(s_0,t^s_0,\dots,s_{k-1}t^s_{k-1})$. We have
$\mu_P([r^s])=\prod_{i<k}2^{-M_i}$. Also,
$\rho^{-1}[[r^s]]=[s]\times[t^s]$ and
$g^{-1}[\rho^{-1}[[r^s]]]=g^{-1}[[s]\times[t^s]]=[s]$, so
$\mu\left(g^{-1}[\rho^{-1}[[r^s]]]\right)=\mu([s])=\prod_{i<k}2^{-M_i}$.
Thus, the claim holds for basic open sets of the form $[r^s]$.
Every clopen subset of $P$ is partitioned by such sets, so the
claim holds for all clopen sets, and hence, as in the proof of
Claim~\ref{rho is measure preserving}, for all Borel sets. This
proves the claim.
 \qed

\m

Now we prove Theorem~\ref{measure analog} in the case that
$B\subseteq 2^\omega$ has outer measure one. It follows from the
usual Lowenheim-Skolem arguments that if we let
$B_\alpha=V^{\poset_\alpha}\cap B$ then there will exist
$\alpha<\omega_2$ such that $B_\alpha\in V^{\poset_\alpha}$ and
$$V^{\poset_\alpha}\models B_\alpha \mbox{ has outer measure one}.$$
Letting $A=\rho^{-1}(B_\alpha)$ (which has outer measure one by
Claim~\ref{rho is measure preserving}) in
\cite[Proposition~2.6]{RS} cited above and using
Claim~\ref{identification of measures}, we have that
$$V^{\poset_{\alpha+1}}\models B_\alpha \mbox{ has $\mu_P$ outer
measure one.}$$ Because ${\mathbb Q}$ is proper, the remainder
$\poset_{\o_2}/\poset_{\a+1}$ of the forcing is isomorphic in
$V^{\poset_{\a+1}}$ to a countable support iteration of ${\mathbb
Q}$ and hence preserves outer measure. It follows that in the
final model $V^{\poset_{\omega_2}}$, $B$ has $\mu_P$ outer measure
one.

Now in the case that $B$ has outer measure less than one, replace
it by $B'=Q+B$ where $Q$ is a countable dense subset of
$2^\omega$. Then $B'$ has outer measure one, and so we know there
exists a measure zero perfect $P$ such that $B'$ has positive
$\mu_P$ outer measure. Hence for some $q\in Q$ we have that $q+B$
has positive $\mu_P$ outer measure.  But then, $B$ has positive
$\mu_{q+P}$ outer measure.

This completes the proof of the theorem.
 \qed

\end{document}